\documentclass[11pt]{amsart}

\usepackage{graphicx}

\newtheorem{theorem}{Theorem}
\newtheorem{corollary}{Corollary}
\newtheorem{lemma}{Lemma}

\def\maxdeg{\operatorname{max-deg}}
\def\maxsupp{\operatorname{max-supp}}
\def\minsupp{\operatorname{min-supp}}
\def\tb{\operatorname{tb}}
\def\sl{\operatorname{sl}}
\def\ssl{\mathfrak{sl}}
\def\maxtb{\overline{\tb}}
\def\maxsl{\overline{\sl}}
\def\P{\mathcal{P}}

\def\Z{\mathbb{Z}}
\def\fig#1{\raisebox{-1.5ex}{\includegraphics[height=4ex]{#1}}}
\def\R{\mathbb{R}}
\def\HKh{\textit{HKh}}
\def\KhB{\textit{KhB}}

\begin{document}

\title{A skein approach to Bennequin type inequalities}
\author{Lenhard Ng}
\address{Mathematics Department, Duke University, Durham, NC 27708}
\email{ng@math.duke.edu}
\urladdr{http://www.math.duke.edu/\~{}ng/}
\thanks{The author is supported by NSF grant DMS-0706777.}

\begin{abstract}
We give a simple unified proof for several disparate bounds on
Thurston--Bennequin number for Legendrian knots and self-linking
number for transverse knots in $\mathbb{R}^3$, and provide a template for
possible future bounds. As an application, we give sufficient
conditions for some of these bounds to be sharp.
\end{abstract}

\maketitle

\section{Introduction}
\label{sec:intro}

\subsection{Main results}

The problem of finding upper bounds for the Thurston--Bennequin and
self-linking numbers of knots has garnered a fair bit of recent attention. Although this originated as a problem in contact geometry, it
now lies more in the realm of knot theory and braid theory, with upper
bounds given by the Seifert and slice genus, the Kauffman and
HOMFLY-PT polynomials, and, more recently, Khovanov homology,
Khovanov--Rozansky homology, and knot Floer homology. These bounds are
sometimes collectively called ``Bennequin type inequalities''.

The original proofs of many Bennequin type inequalities were
remarkably diverse and sometimes somewhat ad hoc. In this paper, we
provide a template which simultaneously proves a number of the
significant Bennequin type inequalities, thus providing a unified
approach to many of these bounds. The proof of the template itself is
a fairly easy induction argument based on the remarkable work of
Rutherford \cite{bib:Ru}.  Our template gives a means to prove future,
yet to be discovered Bennequin type bounds, for example using Khovanov
and Rozansky's proposed categorification of the Kauffman polynomial
\cite{bib:KRKauff}. It also sheds some light on why particular bounds
may be sharp for a certain Legendrian knot while others are not; see
Section~\ref{sec:tree}.

We briefly recall the relevant definitions; see also \cite{bib:Et} or
any number of other references.  A \textit{Legendrian} knot or link in
$\R^3$ with the standard contact structure is a smooth, oriented knot
or link along which $y = dz/dx$ everywhere. It is convenient to
represent Legendrians by their \textit{front projections} to the $xz$
plane, which are (a collection of) oriented closed curves with no
vertical tangencies, whose only singularities are transverse double
points and semicubical cusps. Given a front, one obtains a link
diagram by smoothing out cusps and resolving each double point to a
crossing where the strand with larger slope lies below the strand with
smaller slope; this link diagram, which we will call the
\textit{smoothed front}, is the topological type of the original
Legendrian link. Any topological link has a Legendrian representative.

Given a front $F$, let $c(F)$ denote half of the number of cusps of
$F$, and let $c_\downarrow(F)$ denote the number of cusps of $F$ which
are oriented downwards. Also let $w(F)$ denote the writhe of the
corresponding smoothed front, the number of crossings counted with
the standard signs ($+1$ for
\raisebox{2.5ex}{\includegraphics[height=4ex,angle=270]{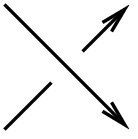}}, $-1$ for
\fig{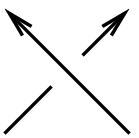}). Define the \textit{Thurston--Bennequin number} and
\textit{self-linking number} of $F$, respectively, by
\begin{align*}
\tb(F) &= w(F) - c(F) \\
\sl(F) &= w(F) - c_\downarrow(F).
\end{align*}
The Thurston--Bennequin and self-linking numbers are invariants of
Legendrian links and comprise the ``classical invariants'' for
Legendrians in $\R^3$. (In the literature, the role of the
self-linking number for Legendrian links is usually played by the rotation
number $\operatorname{r}(F) = \tb(F) - \sl(F)$, and the self-linking
number is reserved for transverse links; our self-linking number for a
Legendrian is the usual self-linking number of its positive transverse
pushoff.)

Within a topological type, $\tb$ and $\sl$ for Legendrian
representatives is always unbounded below; one can decrease $\tb$ and
$\sl$ by adding zigzags to a front. However, in the early 1980's, Bennequin
\cite{bib:Ben}
proved the remarkable fact that $\tb$ and $\sl$ are bounded above for
any link type $L$, by negative the minimal Euler characteristic for a
Seifert surface bounding $L$; for knots $K$,
\[
\tb(F),\sl(F) \leq 2g(K)-1
\]
for any front $F$ representing $K$, where $g(K)$ is the Seifert genus.

Bennequin's inequality has subsequently been improved to a menagerie
of different bounds on $\tb$ and $\sl$. For a topological link $L$,
define the \textit{maximal Thurston--Bennequin number} $\maxtb(L)$
(respectively \textit{maximal self-linking number} $\maxsl(L)$) to be
the maximum
$\tb$ (respectively $\sl$) over all Legendrian realizations of $L$.
Then we have the following Bennequin type inequalities for knots $K$:
\[
\fbox{
$\displaystyle
\begin{aligned}
\maxsl(K) &\leq 2g(K) - 1 && \text{(Bennequin bound \cite{bib:Ben})}
\\
\maxsl(K) &\leq 2g_4(K) - 1 && \text{(slice-Bennequin bound \cite{bib:Rud2})}
\\
\maxsl(K) &\leq 2\tau(K) - 1 && \text{($\tau$ bound
  \cite{bib:Plam1})}
\\
\maxsl(K) &\leq s(K) - 1 && \text{($s$ bound \cite{bib:Plam2,bib:Shu})}
\\
\maxsl(K) &\leq -\maxdeg_a P(K)(a,z)-1 && \text{(HOMFLY-PT bound
  \cite{bib:FW,bib:Mor})}
\\
\maxsl(K) &\leq -\maxsupp_a \P^{a,q,t}(K)-1 && \text{(HOMFLY-PT homology
  bound \cite{bib:Wu})}
\\
\maxtb(K) &\leq -\maxdeg_a F(K)(a,z)-1 && \text{(Kauffman bound
  \cite{bib:Rud1})}
\\
\maxtb(K) &\leq -\minsupp_{q-t} \HKh^{q,t}(K) && \text{(Khovanov bound
  \cite{bib:NgKho})}
\end{aligned}
$}
\]

For notation and conventions, see Section~\ref{ssec:notation}.
Several remarks are in order.

\begin{itemize}

\item
Most of these inequalities have obvious generalizations to links; in
particular, the last four translate unchanged to bounds for links.

\item $\maxtb(L)
\leq \maxsl(L)$ always: rotating any front $F$
$180^\circ$ produces another front $F'$ of the same topological type
with $2\tb(F) = 2\tb(F') = \sl(F) + \sl(F')$. It follows that any
upper bound for $\maxsl$ is also an upper bound for $\maxtb$. However,
the Kauffman and Khovanov $\maxtb$
bounds above do not extend to bounds on $\maxsl$.

\item Some Bennequin type inequalities imply others. The $\tau$ and
$s$ bounds (and presumably the HOMFLY-PT homology bound) imply
slice-Bennequin, which in turn implies Bennequin; the HOMFLY-PT
homology bound also implies the HOMFLY-PT (polynomial) bound. On the
other hand, many pairs of the inequalities are incommensurable,
notably the Kauffman and Khovanov bounds \cite{bib:NgKho} (see also
\cite{bib:Fer}).

\item The above inequalities (in particular, the Kauffman, Khovanov,
  and HOMFLY-PT bounds) suffice to calculate $\maxtb$ and $\maxsl$ for
  all but a handful of knots with $11$ or fewer crossings \cite{bib:Ngarc}.

\item The Kauffman and HOMFLY-PT bounds have been given many proofs in
  the literature, involving state models, plane curves, the Jaeger
  formula, etc.
  See \cite{bib:CG,bib:Fer,bib:FT,bib:Tab,bib:Tan} for additional proofs
  of the Kauffman bound, and \cite{bib:CGM,bib:FT,bib:Tab} for HOMFLY-PT.

\item This is not a complete list of known Bennequin type
inequalities. In particular, Wu \cite{bib:Wu,bib:Wu3} has derived
bounds on $\maxsl$ from Khovanov--Rozansky $\ssl_n$ homology.

\end{itemize}

Our main results give general criteria for a link invariant to provide
an upper bound for $\maxtb$ or $\maxsl$. These criteria
are satisfied for many of the known bounds.

\begin{theorem}
Suppose \label{thm:sl} that $i(L)$ is a $\Z$-valued invariant of oriented
links such that:
\begin{enumerate}
\item
$i( \overbrace{\bigcirc \bigcirc \dots \bigcirc}^n) \leq n$;
\label{cond:sl1}
\item
we have \label{cond:sl2}
\[
i(\fig{skein1.eps})+1 \leq \max\left( i(\fig{skein2.eps})-1,
i(\fig{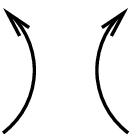})\right)
\]
and
\[
i(\fig{skein2.eps})-1 \leq \max\left( i(\fig{skein1.eps})+1,
i(\fig{skein3.eps})\right).
\]
\end{enumerate}
Then
\[
\maxsl(L) \leq -i(L).
\]
\end{theorem}

\begin{corollary}
The HOMFLY-PT and
HOMFLY-PT homology bounds on $\maxsl$ hold for oriented links.
\label{cor:HH}
\end{corollary}

\begin{theorem}
Suppose that $\tilde{i}(D)$ is a $\Z$-valued invariant of unoriented
link diagrams such that:
\begin{enumerate}
\item \label{cond:tb1}
$\tilde{i}(D)$ is invariant under Reidemeister moves II and III;
\item \label{cond:tb2}
$\tilde{i}(\fig{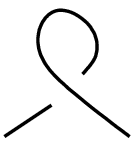}) = \tilde{i}(\fig{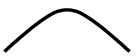}) + 1$;
\item \label{cond:tb3}
$\tilde{i}( \overbrace{\bigcirc \bigcirc \dots \bigcirc}^n) \leq n$;
\item
\label{cond:tb4}
we have
\[
\tilde{i}(\fig{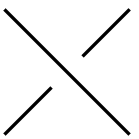}) \leq \max\left(
\tilde{i}(\fig{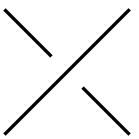}),
\tilde{i}(\fig{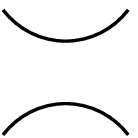}),
\tilde{i}(\fig{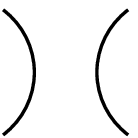}) \right).
\]
\end{enumerate}
Then \label{thm:tb}
\[
\maxtb(L) \leq -i(L),
\]
where $i(L)$ is the link invariant defined by $i(D) = \tilde{i}(D) -
w(D)$.
\end{theorem}

\begin{corollary}
The Kauffman and Khovanov bounds on $\maxtb$ hold for oriented links.
\label{cor:KK}
\end{corollary}

We speculate that when $\tau$
and $s$ are extended to oriented links, $i=1-2\tau$ and $i=1-s$ should
each satisfy the conditions of Theorem~\ref{thm:sl}. This would demonstrate
that the $\tau$ and $s$ bounds can be proven using our template as
well.

It also seems likely that Khovanov--Rozansky's proposed
categorification of the Kauffman polynomial \cite{bib:KRKauff} would
satisfy a skein relation which would allow one to apply
Theorem~\ref{thm:tb}. This would give an upper bound on $\maxtb$ from
Kauffman homology strengthening the Kauffman (polynomial) bound, just
as the HOMFLY-PT homology bound on $\maxsl$ strengthens the HOMFLY-PT
(polynomial) bound.

As mentioned earlier, one benefit of our results is a better
understanding of when particular Bennequin type inequalities are
sharp. Rutherford \cite{bib:Ru} has demonstrated a necessary and
sufficient condition for the Kauffman bound to be sharp, in terms of
certain decompositions of fronts known as rulings. It would be
nice to have similar characterizations for sharpness for, say,
the HOMFLY-PT bound and the Khovanov bound. It seems that such
characterizations should now be within reach, but for now we present
some sufficient conditions for these bounds to be sharp; see
Section~\ref{sec:tree}.

Here is a rundown of the rest of the paper.
In Section~\ref{ssec:notation}, we summarize the notation used in our
presentation of the Bennequin type inequalities. The proofs of
Theorems~\ref{thm:sl} and \ref{thm:tb} and their rather easy
consequences, Corollaries~\ref{cor:HH} and \ref{cor:KK}, are given in
Section~\ref{sec:proofs}. In Section~\ref{sec:tree}, we use the
inductive proofs of our main results to construct trees which
decompose any Legendrian link into simpler links, and use these trees
to study sharpness of some Bennequin type inequalities.

\subsection{Notation}
\label{ssec:notation}

Here we collect the definitions used in the Bennequin type
inequalities mentioned above, including the particular conventions we
use. (These conventions coincide with those from \texttt{KnotTheory}
\cite{bib:BN} wherever applicable.)

\begin{itemize}

\item $\maxdeg_a$ is the maximum degree in $a$; $\maxsupp_a$ is the
  maximum $a$ degree in which the homology is supported;
  $\minsupp_{q-t}$ is the minimum value for $q-t$ over all bidegrees
  $(q,t)$ in which the homology is supported. 
\item
$g_4(K)$ is the slice genus of $K$.
\item
$\tau(K)$ is the concordance invariant from knot
Floer homology \cite{bib:OSz}, normalized so that $\tau = 1$ for
right-handed trefoil.
\item
$s(K)$ is Rasmussen's concordance invariant from
Khovanov homology \cite{bib:Ras}.
\item
$P(K)(a,z)$ is the HOMFLY-PT polynomial of $K$,
normalized so that $P=1$ for the unknot and
\[
aP(\fig{skein1.eps}) - a^{-1} P(\fig{skein2.eps}) = z P(\fig{skein3.eps}).
\]
\item
$\P^{a,q,t}(K)$ is the (reduced) Khovanov--Rozansky
HOMFLY-PT homology \cite{bib:KhR2}
categorifying $P(K)$, normalized so that
\[
\sum_{i,j,k} (-1)^{(k-i)/2} a^i q^j \dim \P^{i,j,k}(K) = P(K)(a,q-q^{-1}).
\]
Note that our $\P^{a,q,t}(K)$ is Rasmussen's
$\overline{H}^{q,-a,t}(K)$ from \cite{bib:RasKR}. (For links, we need
to use the ``totally reduced'' version $\overline{\overline{H}}$
of HOMFLY-PT homology in place
of reduced HOMFLY-PT homology $\overline{H}$; see \cite{bib:RasKR}.)
\item
$F(K)(a,z)$ is the Kauffman polynomial of $K$,
normalized so that for a diagram $D$ representing $K$,
$F(K)(a,z)=a^{-w(D)} \tilde{F}(D)(a,z)$, where
$\tilde{F}$ is the framed Kauffman polynomial, the regular-isotopy invariant of
unoriented link diagrams defined by $\tilde{F}(\bigcirc)
= 1$, $\tilde{F}(\fig{pluskink.eps}) = a\tilde{F}(\fig{plusnokink.eps})$, and
\[
\tilde{F}(\fig{unskein1.eps}) + \tilde{F}(\fig{unskein2.eps}) =
z\left(\tilde{F}(\fig{unskein3.eps}) +
\tilde{F}(\fig{unskein4.eps})\right).
\]
\item
$\HKh^{q,t}(K)$ is ($\ssl_2$) Khovanov homology,
normalized so that
\[
\sum_{i,j} (-1)^j q^i \dim \HKh^{i,j}(K) = (q+q^{-1}) V_K(q^2),
\]
where $V_K$ is the Jones polynomial.
\end{itemize}

\section{Proofs}
\label{sec:proofs}

Theorems~\ref{thm:sl} and \ref{thm:tb} have essentially the same
proof. We establish Theorem~\ref{thm:tb} first, and then prove
Theorem~\ref{thm:sl} and Corollaries~\ref{cor:HH} and~\ref{cor:KK}.

\begin{proof}[Proof of Theorem~\ref{thm:tb}]

View $\tilde{i}$ as a map on fronts by applying $\tilde{i}$ to the
smoothed version of any front.
Let $L$ be an oriented link and let $F$ be a Legendrian front of type
$L$. We wish to show that $\tb(F) \leq -i(L)$, or equivalently, that
$c(F) - \tilde{i}(F) \geq 0$. Note that $c(F) - \tilde{i}(F)
= - \tb(F)-i(F)$ is invariant under Legendrian isotopy.

The idea, which is essentially due to Rutherford \cite{bib:Ru}, is to
use skein moves to replace $F$ by simpler fronts in such a way that
$c-\tilde{i}$ does not decrease, and then to induct. The four fronts
\fig{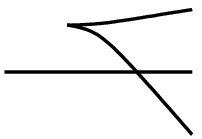}, \fig{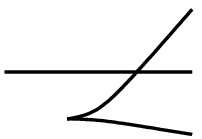}, \fig{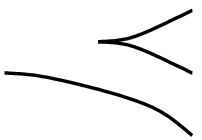}, and
\fig{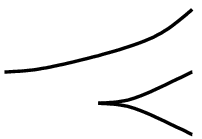} are topologically
\fig{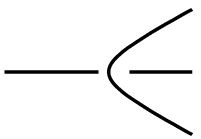}, \fig{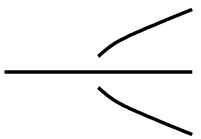},
\fig{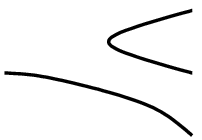}, and \fig{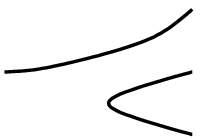}, respectively, and
are thus related by the four-term unoriented skein relation.

Suppose that $F$ contains a Legendrian tangle \fig{Legskein1.eps}.
Successively replace
\begin{equation}
\fig{Legskein1.eps} ~~~~~\longrightarrow
~~~~~\fig{Legskein2.eps}~~~~~,~~~~~
\fig{Legskein3.eps}~~~~~,~~~~~
\fig{Legskein4.eps}
\label{eq:skeinreplace1}
\end{equation}
in $F$, to obtain three new fronts, and
suppose that $c-\tilde{i} \geq 0$ for each of those fronts. Since $c$
is the same for all four fronts, assumption (\ref{cond:tb4})
in the statement of Theorem~\ref{thm:tb}
then implies that $c-\tilde{i} \geq 0$ for $F$ as well. Similarly, if
$F$ contains \fig{Legskein2.eps}, and the three fronts obtained from
$F$ by
\begin{equation}
\fig{Legskein2.eps} ~~~~~\longrightarrow ~~~~~\fig{Legskein1.eps}~~~~~,~~~~~
\fig{Legskein3.eps}~~~~~, ~~~~~\fig{Legskein4.eps}
\label{eq:skeinreplace2}
\end{equation}
all satisfy $c-\tilde{i} \geq 0$, then $c-\tilde{i} \geq 0$ for $F$ as well.

To prove that $c(F)-\tilde{i}(F) \geq 0$ for all $F$,
we induct on the \textit{singularity number} $s(F)$ of $F$, defined as the
total number of singularities (crossings and cusps) of $F$.
If $s(F)=2$, then $F$ is the standard
Legendrian unknot and $c(F) = 1$, $\tilde{i}(F) \leq 1$. Now consider a
general front $F$. Suppose that $F$ contains a tangle of the form
\fig{Legskein1.eps} or \fig{Legskein2.eps}. If we replace this tangle
successively by three tangles according to (\ref{eq:skeinreplace1}) or
(\ref{eq:skeinreplace2}), then the last two of the resulting fronts
have lower $s$ than $F$ and are covered by the
induction assumption, while the first has the same $s$
as $F$.

The strategy is now to apply ``skein crossing changes''
$\fig{Legskein1.eps} \leftrightarrow \fig{Legskein2.eps}$
to obtain a simpler
front. To do this, we perform a second induction, this time on a
modified singularity number $s'(F)$, defined as the number of
singularities to the right of the rightmost left cusp of $F$.
Since Legendrian isotopy, Legendrian destabilization (the removal of
a zigzag), and the removal of trivial unknots
do not increase $c-\tilde{i}$, the Theorem follows by induction from
the following result.

\begin{lemma}[Rutherford \cite{bib:Ru}, Lemma 3.3]
Via skein crossing changes, Legendrian isotopy, Legendrian
destabilization, and the removal of trivial unknots,
we can turn $F$ into a front which either has lower
$s$, or the same $s$ and lower $s'$.
\label{lem}
\end{lemma}

\begin{figure}
\centerline{
\includegraphics[height=7in]{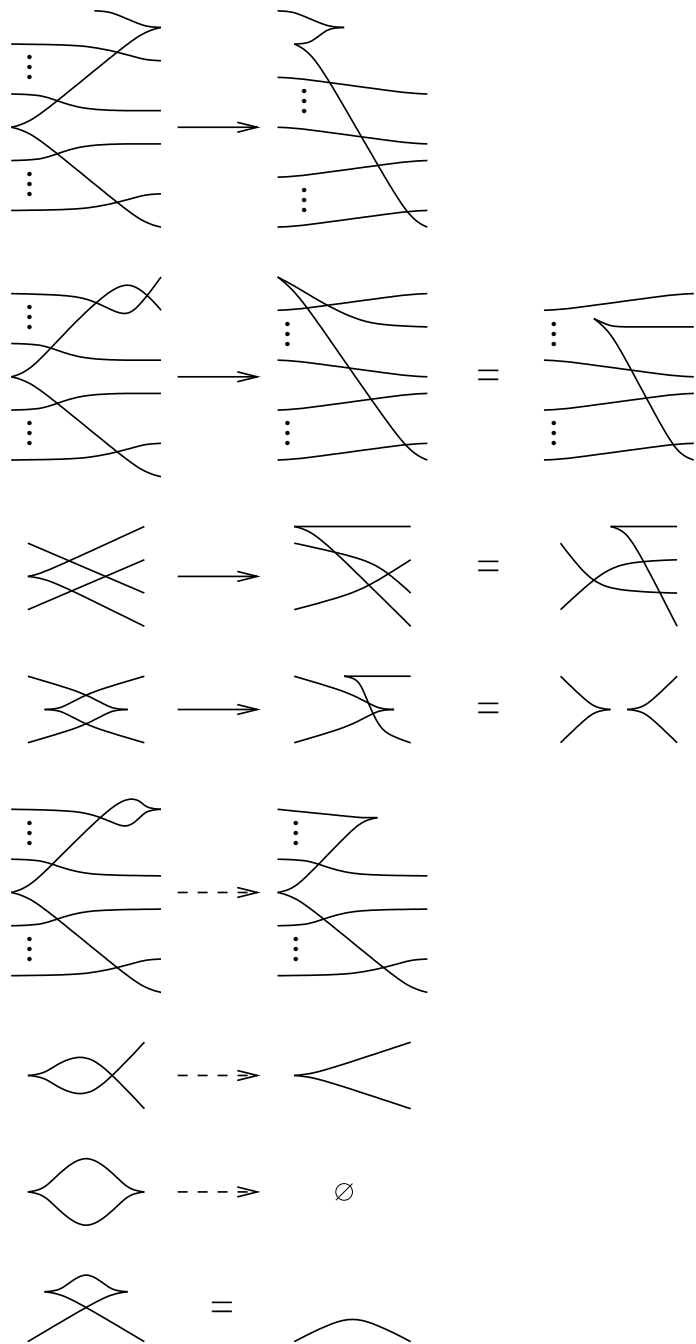}
}
\caption{
Inductively simplifying fronts. Solid arrows denote skein crossing
changes; dashed arrows denote destabilizations or deletions of trivial
unknots;
equalities denote Legendrian isotopies. For each of the initial
fronts (left column),
the rightmost left cusp in the front is the unique left cusp depicted.
}
\label{fig:maininduct}
\end{figure}

For completeness, we sketch here the proof of the lemma.
Consider the portion of $F$ immediately to the right
of the rightmost left cusp of $F$. By using Legendrian Reidemeister
moves II and III if necessary, we can assume that this portion of $F$
has one of the forms shown on the left hand side of
Figure~\ref{fig:maininduct}. In each case, the use of skein crossing
changes, Legendrian isotopy, Legendrian destabilization, and/or the
removal of trivial unknots yields a simpler front (one with lower $s$,
or the same $s$ and lower $s'$). The lemma, and Theorem~\ref{thm:tb}, follows.
\end{proof}

\begin{proof}[Proof of Theorem~\ref{thm:sl}]

This is a minor modification of the proof of
Theorem~\ref{thm:tb}. Define an invariant $\tilde{i}$ of oriented link
diagrams by $\tilde{i}(D) = i(D) + w(D)$; then
\[
\tilde{i}(\fig{skein1.eps}) \leq \max\left( \tilde{i}(\fig{skein2.eps}),
\tilde{i}(\fig{skein3.eps}) \right)
\]
and
\[
\tilde{i}(\fig{skein2.eps}) \leq \max\left( \tilde{i}(\fig{skein1.eps}),
\tilde{i}(\fig{skein3.eps}) \right).
\]
For (oriented) Legendrian fronts $F$, we wish to show that $\sl(F) \leq -i(F)$,
or equivalently, that $c_\downarrow(F) - \tilde{i}(F) \geq 0$.

As before, we use skein moves to induct on the singularity number of
$F$. Note that $c_\downarrow(F) - \tilde{i}(F)$ is invariant under
Legendrian isotopy and nonincreasing under Legendrian destabilization.
If $F$ contains a tangle \fig{Legskein1.eps} (respectively
\fig{Legskein2.eps}), then we can successively
replace it by \fig{Legskein2.eps} (respectively \fig{Legskein1.eps})
and whichever of \fig{Legskein3.eps}
and \fig{Legskein4.eps} inherits an orientation from $F$, to obtain
two new fronts. If $c_\downarrow-\tilde{i} \geq 0$ for these two
fronts, then $c_\downarrow-\tilde{i} \geq 0$ for $F$ as well. We now
apply Lemma~\ref{lem} as before.
\end{proof}

\begin{proof}[Proof of Corollary~\ref{cor:HH}]
Define $i(L) = \maxdeg_a P(L)(a,z) + 1$. By the skein relation and
normalization for the HOMFLY-PT polynomial, the conditions in
Theorem~\ref{thm:sl} hold, and Theorem~\ref{thm:sl} then gives the
HOMFLY-PT bound.

For the HOMFLY-PT homology bound, define $i(L) = \maxsupp_a
\P^{a,q,t}(L) + 1$.
The skein relation for the HOMFLY-PT
polynomial (see \cite{bib:RasKR}) categorifies to an exact triangle relating
$\P(\fig{skein1.eps})$, $\P(\fig{skein2.eps})$, and
$\P(\fig{skein3.eps})$, and this exact triangle yields condition
(\ref{cond:sl2}) in the statement of Theorem~\ref{thm:sl}. The
normalization condition (\ref{cond:sl1}) is easy to check, and thus
Theorem~\ref{thm:sl} yields the HOMFLY-PT homology bound.
\end{proof}

\begin{proof}[Proof of Corollary~\ref{cor:KK}]
For the Kauffman bound, define $\tilde{i}(D) = \maxdeg_a \tilde{F}(D)
+ 1$ for unoriented link diagrams $D$; by the skein relation for
$\tilde{F}$, the conditions for Theorem~\ref{thm:tb} are satisfied,
and the Kauffman bound follows.

For the Khovanov bound, collapse the $q,t$-bigrading on Khovanov homology to
a single grading $* = q-t$. One can define framed Khovanov homology
$\widetilde{\HKh}^*(D)$ for unoriented link diagrams $D$ such that
$\HKh^*(L) = \widetilde{\HKh}^{*-w(D)}(D)$ if $D$ is of link type $L$;
indeed, the complex for Khovanov homology is first defined this
way. The exact triangle
in Khovanov homology (in this context, see, e.g., \cite{bib:NgKho}) is
given by
\[
\cdots \longrightarrow \widetilde{\HKh}^*(\fig{unskein4.eps})
\longrightarrow
\widetilde{\HKh}^*(\fig{unskein1.eps}) \longrightarrow
\widetilde{\HKh}^*(\fig{unskein3.eps}) \longrightarrow
\widetilde{\HKh}^{*-1}(\fig{unskein4.eps}) \longrightarrow \cdots.
\]

If we now define $\tilde{i}(D) = -\minsupp_* \widetilde{\HKh}^*(D)$,
then the exact triangle implies that
\[
\tilde{i}(\fig{unskein1.eps}) \leq \max \left(
\tilde{i}(\fig{unskein3.eps}),
\tilde{i}(\fig{unskein4.eps}) \right).
\]
In particular, condition (\ref{cond:tb4}) in Theorem~\ref{thm:tb}
holds, and the Khovanov bound follows.
\end{proof}


\section{Skein trees}
\label{sec:tree}

The nature of the proofs of Theorems~\ref{thm:sl} and~\ref{thm:tb}
allows us to give necessary conditions and sufficient conditions for
various Bennequin type inequalities to be sharp, and to compare these
inequalities with each other. We can decompose any Legendrian knot via
a skein tree, much as one would do to calculate knot polynomials using
skein relations, and the skein tree can often tell us whether one
bound or another is sharp.

\begin{figure}
\centerline{
\includegraphics[width=5.5in]{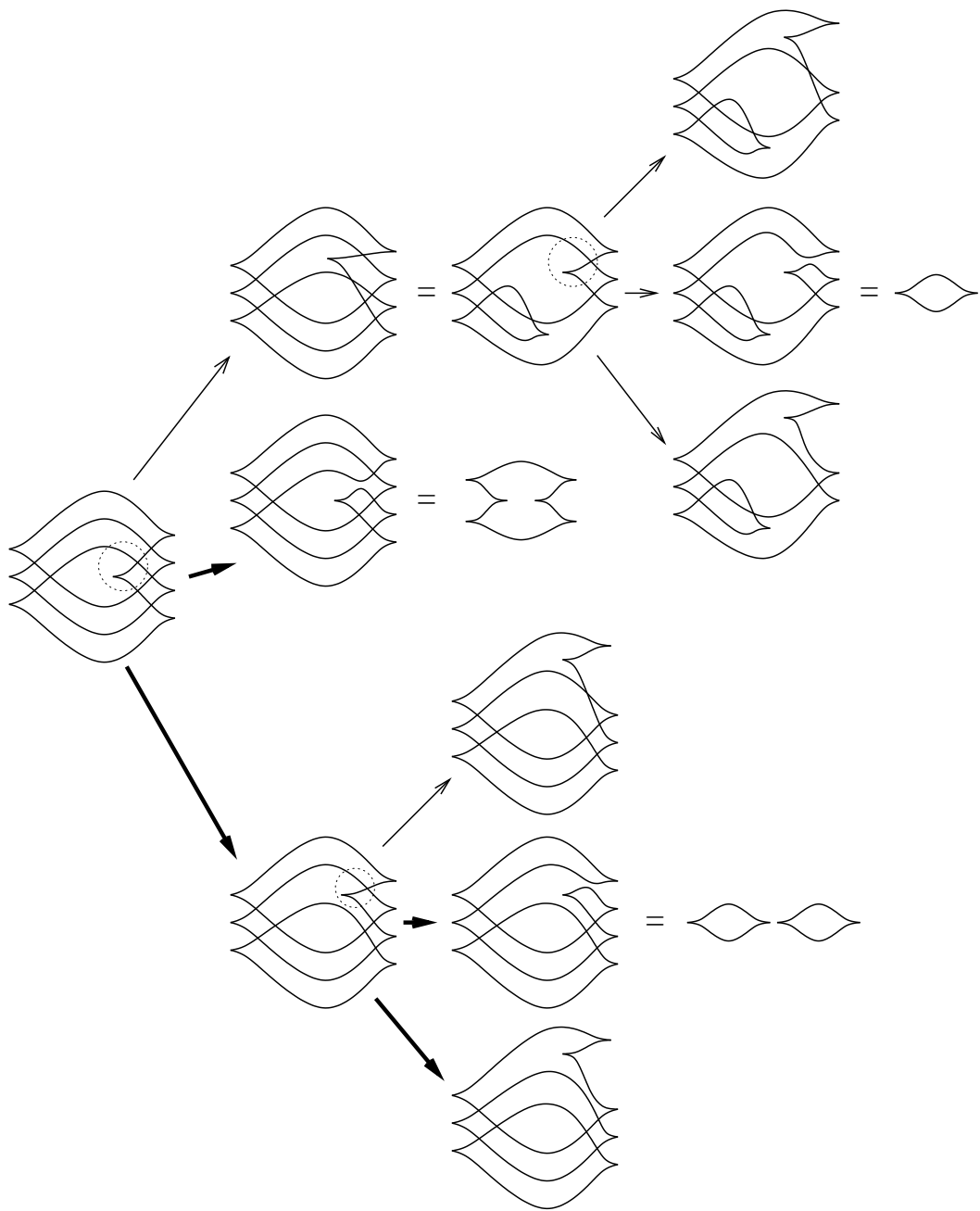}
}
\caption{
Unoriented skein tree for a Legendrian $(3,-4)$ torus knot. The
Jones skein tree is the subtree given by the bold arrows. The circled tangles are where the tangle replacements are made.
}
\label{fig:torusskeintree}
\end{figure}

Starting with an unoriented Legendrian front, construct the
\textit{unoriented skein tree} by following Rutherford's strategy
described in the proof of Theorem~\ref{thm:tb}:
\begin{itemize}
\item
at each step, do a tangle replacement
\[
\fig{Legskein1.eps} ~~~~~\longrightarrow
~~~~~\fig{Legskein2.eps}~~~~~,~~~~~
\fig{Legskein3.eps}~~~~~,~~~~~
\fig{Legskein4.eps}
\]
or
\[
\fig{Legskein2.eps} ~~~~~\longrightarrow ~~~~~\fig{Legskein1.eps}~~~~~,~~~~~
\fig{Legskein3.eps}~~~~~, ~~~~~\fig{Legskein4.eps}
\]
to obtain three new fronts;
\item
simplify the results by Legendrian isotopy, and repeat;
\item
stop when the result is either a stabilization (isotopic to a front with a zigzag) or a standard Legendrian
unlink (the disjoint union of $\tb=-1$ Legendrian unknots).
\end{itemize}
An example is given in Figure~\ref{fig:torusskeintree}.

One can easily use an unoriented skein tree for $F$ to calculate the
coefficient of $a^{-\tb(F)-1}$ in the Kauffman polynomial for $F$ with any orientation
(this coefficient is nonzero if and only the Kauffman bound is sharp):
for terminal leaves in the tree, the coefficient is $1$ at a
standard Legendrian unlink and $0$ at a stabilized front; use the skein
relation for the framed Kauffman polynomial to backwards-construct
the coefficient along the tree:
\[
\fig{Legskein1.eps} + \fig{Legskein2.eps} = z\left(
  \fig{Legskein3.eps} + \fig{Legskein4.eps} \right).
\]
This is simply a restatement of a result of Rutherford \cite{bib:Ru}.

We now see a heuristic reason for why the Kauffman bound sometimes
fails even for Legendrian knots which maximize $\tb$. Consider for example
the Legendrian $(3,-4)$ torus knot $F$ shown in
Figure~\ref{fig:torusskeintree}. Two of the terminal leaves of the
unoriented skein tree are standard Legendrian unlinks. Either leaf by
itself gives a contribution to the Kauffman polynomial of $T(3,-4)$
which would imply that the Kauffman bound on $F$ is sharp; but the two
contributions cancel, and the maximum framing degree of the Kauffman
polynomial of $T(3,-4)$ is one less than necessary for sharpness.

%

Next we examine the sharpness of the Khovanov bound on $\tb$. For a
front $F$, define
\[
\KhB(F) = -\tb(F)+\minsupp_{q-t} \HKh^{q,t}(F)
\]
(this is $c-\tilde{i}$ from the proof of Theorem~\ref{thm:tb}); then
$\KhB(F) \geq 0$ with equality if and only if the Khovanov bound is
sharp for $F$.

Define the \textit{Jones skein tree} to be the subtree of the
unoriented skein tree consisting only of tangle replacements
\[
\fig{Legskein1.eps} ~~~~~\longrightarrow
~~~~~\fig{Legskein3.eps}~~~~~,~~~~~
\fig{Legskein4.eps}
\]
and
\[
\fig{Legskein2.eps} ~~~~~\longrightarrow ~~~~~
\fig{Legskein4.eps}~~~~~, ~~~~~\fig{Legskein3.eps}.
\]
(In each case, the two replacements are the $0$- and $1$-resolution, respectively.) Since this only counts $0$- and $1$-resolutions and not crossing changes, this is the same tree used to calculate the Jones polynomial for a knot.

At each stage in the Jones skein tree, a front $F$ is connected to its
$0$-resolution $F_0$ and its $1$-resolution $F_1$. The skein exact
sequence for Khovanov homology implies the following.

\begin{lemma}
If $\KhB(F_0) = 0$, or $\KhB(F_1) = 0$ and $\KhB(F_0) \geq 2$, then
$\KhB(F) = 0$.
\label{lem:Khsharp}
\end{lemma}

We now have the following sufficient condition for the Khovanov bound
to be sharp.

\begin{theorem}
Let $F$ be a front.
Circle particular fronts in the Jones skein tree for $F$ as follows:
circle all terminal leaves which are standard Legendrian unknots; then
work backwards, circling a front $F'$ if either
\begin{itemize}
\item
$(F')_0$ is circled or
\item
$(F')_1$ is circled and $(F')_0$ is isotopic to a front stabilized at
least twice (i.e., with two zigzags).
\end{itemize}
If $F$ is circled by this process, then the Khovanov bound is sharp
for $F$.
\label{thm:Khsharp}
\end{theorem}

\begin{proof}
Since $\KhB \geq 2$ for a front stabilized at least twice,
Lemma~\ref{lem:Khsharp} implies that all circled fronts satisfy $\KhB
= 0$.
\end{proof}

\begin{figure}
\centerline{
\includegraphics[height=2.5in]{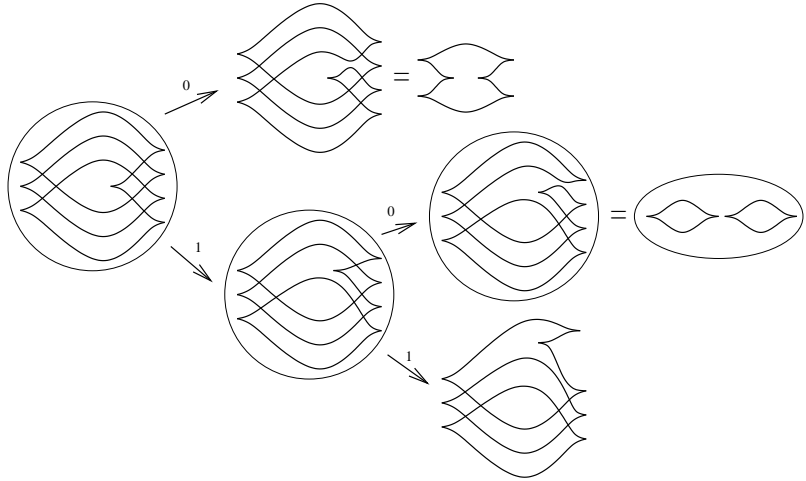}
}
\caption{The Jones skein tree for the $(3,-4)$ torus knot. The numbers
  refer to $0$- or $1$-resolutions; the circled fronts are Khovanov sharp.
}
\label{fig:torusKhtree}
\end{figure}

As an illustration of Theorem~\ref{thm:Khsharp}, the Khovanov bound is
sharp for the Legendrian $(3,-4)$ torus knot shown in
Figure~\ref{fig:torusskeintree}; see
Figure~\ref{fig:torusKhtree}. Comparing with
Figure~\ref{fig:torusskeintree} gives some indication of why the
Khovanov bound is sharp here but the Kauffman bound is not: two
terminal leaves of the unoriented skein tree are standard unlinks and
their contributions to the Kauffman polynomial cancel, while only one
of these terminal leaves is counted in the Jones skein tree and it
makes a nonvanishing contribution to Khovanov homology.

We remark that Theorem~\ref{thm:Khsharp} implies, but
is generally much stronger than, the sufficient condition for Khovanov
sharpness given in \cite{bib:NgKho}. Recall from \cite{bib:NgKho} that
the $0$-resolution of a front, obtained by replacing each double point
by its $0$-resolution, is \textit{admissible} if each component of the
$0$-resolution is a standard Legendrian unknot, and no component
contains both pieces of any resolved double point.

\begin{corollary}[{\cite[Proposition 7]{bib:NgKho}}]
The Khovanov bound is sharp for any front with admissible $0$-resolution.
\end{corollary}

\begin{proof}
Suppose that $F$ is a front with admissible $0$-resolution; apply the
procedure from Theorem~\ref{thm:Khsharp}.
It is easy to check that in the Jones skein tree for $F$, $F$ and all
of its iterated $0$-resolutions are circled.
\end{proof}

\noindent
We do not know how the condition of Theorem~\ref{thm:Khsharp} compares
to the sufficient condition for Khovanov sharpness given by Wu \cite{bib:Wu2}.

\begin{figure}
\centerline{
\includegraphics[width=4.5in]{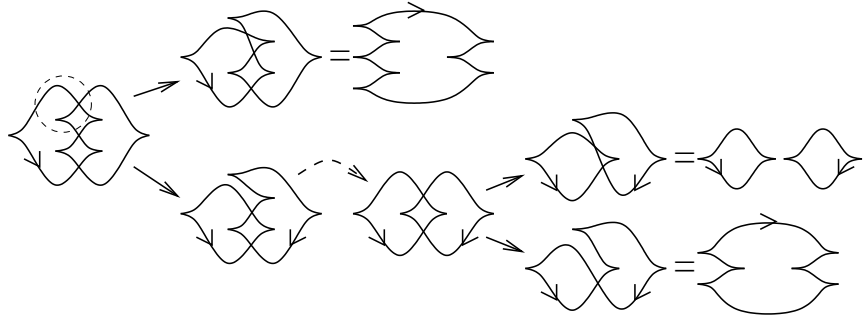}
}
\caption{
Oriented skein tree for a Legendrian trefoil.
Solid arrows represent the skein relation; the dashed arrow is a
negative destabilization.
}
\label{fig:trefoilskeintree}
\end{figure}

One can similarly construct an \textit{oriented skein tree} for any
oriented front, at each step replacing a front by the two fronts
related to it by the oriented skein relation. Rather than stopping at
all stabilized fronts, as for the unoriented skein tree, we stop only at
fronts which are positive stabilizations (i.e., isotopic to a
front with a downward zigzag \fig{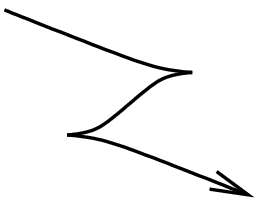}). If we encounter a negative
stabilization (i.e., a front isotopic to one with an upward zigzag
\fig{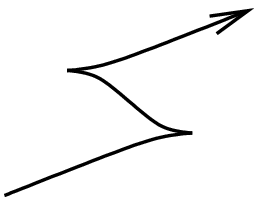}), we eliminate the upward zigzag (this does not change
$\sl$) and proceed. All terminal leaves of the oriented skein
tree are either positive stabilizations or standard Legendrian
unlinks. See Figure~\ref{fig:trefoilskeintree}.

As for unoriented skein trees and the Kauffman polynomial, we can use
an oriented skein tree for a front $F$ to calculate the coefficient of
$a^{-\sl(F)-1}$ in the HOMFLY-PT polynomial for $F$ (this coefficient
is nonzero if and only if the HOMFLY-PT bound is sharp). Standard
Legendrian unlinks have coefficient $1$; positive stabilizations have
coefficient $0$; the coefficient is preserved under negative
stabilizations; we can backwards construct the coefficient along the
tree using the skein relation, e.g.,
\[
\fig{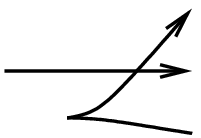} - \fig{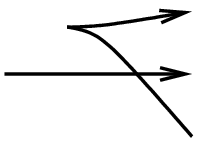} = z\, \fig{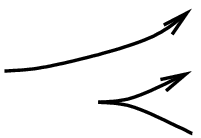}
\hspace{3ex} \text{and} \hspace{3ex}
\fig{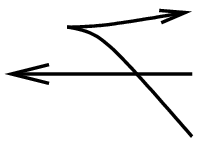} - \fig{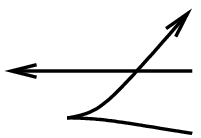} = z \, \fig{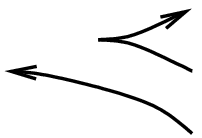}.
\]
This again is very similar to a result of Rutherford \cite{bib:Ru}.

It is sometimes easy to tell by inspection of an oriented skein tree
whether the HOMFLY-PT bound is sharp. For example:

\begin{theorem}
If an odd number of terminal leaves of the oriented skein tree of $F$
are standard Legendrian unlinks, then the HOMFLY-PT polynomial bound
on $\sl$ is sharp for $F$.
\end{theorem}

\noindent
Obviously this sufficient condition is rather weak, but it does for
instance imply that the Legendrian trefoil in
Figure~\ref{fig:trefoilskeintree} maximizes $\sl$.

\begin{figure}
\centerline{
\includegraphics[height=3in]{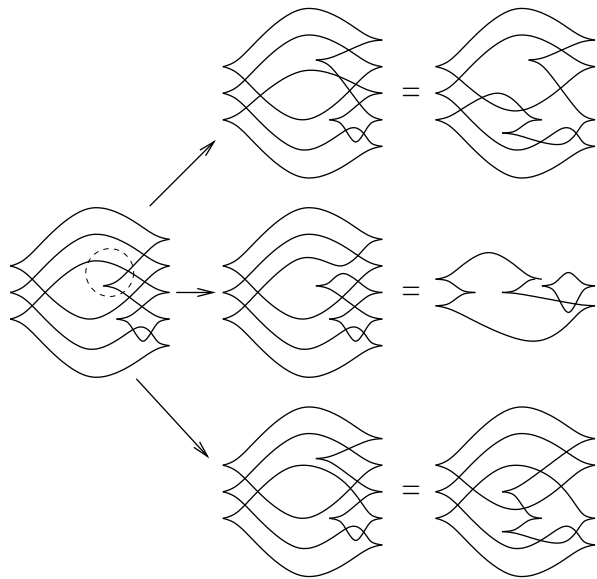}
}
\caption{
Unoriented skein tree for a Legendrian $m(10_{132})$ knot.
}
\label{fig:10132}
\end{figure}

Skein trees can also show the limitations of Theorems~\ref{thm:sl}
and~\ref{thm:tb}. Consider the unoriented skein tree for the
$m(10_{132})$ knot $F$ shown in Figure~\ref{fig:10132}. Each of the
terminal leaves of the tree is a stabilization. Now suppose that
$\tilde{i}(D)$ is any invariant satisfying the conditions of
Theorem~\ref{thm:tb}. Then by Theorem~\ref{thm:tb},
$c-\tilde{i} \geq 0$ for all fronts; since each of the fronts on the
right hand side of Figure~\ref{fig:10132} is a stabilization,
$c-\tilde{i} \geq 1$ for these. Condition (\ref{cond:tb4}) from
Theorem~\ref{thm:tb} implies that $c(F)-\tilde{i}(F) \geq 1$ as well,
and so $i(m(10_{132})) \leq -\tb(F)-1 = 0$.

It follows that the best possible bound given by Theorem~\ref{thm:tb}
is $\maxtb(m(10_{132})) \leq 0$. We however know from
\cite{bib:Ngarc} that $\maxtb(m(10_{132})) = -1$; hence the template
of Theorem~\ref{thm:tb} can never give a sharp bound for
$\maxtb(m(10_{132}))$. A similar argument shows that the template of
Theorem~\ref{thm:sl} can never give a sharp bound for
$\maxsl(m(10_{132}))$ ($=-1$ by \cite{bib:Ngarc}).

\section*{Acknowledgments}

I would like to thank Matt Hedden and Jake Rasmussen for helpful
conversations, and 
Princeton University for its hospitality during the course of this
work.



\begin{thebibliography}{10}

\bibitem{bib:BN}
D.\ Bar-Natan, The Mathematica package
\texttt{KnotTheory\`{}\,}, available at The Knot Atlas,
\verb+http://katlas.math.toronto.edu/wiki/+.

\bibitem{bib:Ben}
D.\ Bennequin, Entrelacements et \'equations de Pfaff,
\textit{Ast\'erisque} \textbf{107--108} (1983), 87--161.

\bibitem{bib:CG}
S.\ Chmutov and V.\ Goryunov, Polynomial invariants of Legendrian
links and their fronts, \textit{KNOTS '96 (Tokyo)},  239--256, World
Sci. Publ., River Edge, NJ, 1997.

\bibitem{bib:CGM}
S.\ Chmutov, V.\ Goryunov, and M.\ Murakami, Regular Legendrian knots
and the HOMFLY polynomial of immersed plane curves, \textit{Math.\
  Ann.} \textbf{317} (2000),  no.\ 3, 389--413.

\bibitem{bib:Et}
J.\ B.\ Etnyre, Legendrian and transversal knots, in \textit{The
Handbook of Knot Theory} (Elsevier, Amsterdam, 2005), 105--185;
\texttt{math/0306256}.

\bibitem{bib:Fer}
E.\ Ferrand, On Legendrian knots and polynomial invariants,
\textit{Proc.\ Amer.\ Math.\ Soc.} \textbf{130} (2002), no.\ 4,
1169--1176; \texttt{math/0002250}.

\bibitem{bib:FW}
J.\ Franks and R.\ F.\ Williams, Braids and the Jones polynomial,
\textit{Trans.\ Amer.\ Math.\ Soc.} \textbf{303} (1987), no.\ 1, 97--108.

\bibitem{bib:FT}
D.\ Fuchs and S.\ Tabachnikov, Invariants of Legendrian and transverse
knots in the standard contact space, \textit{Topology} \textbf{36}
(1997),  no.\ 5, 1025--1053.


\bibitem{bib:Kh}
M.\ Khovanov, A categorification of the Jones
polynomial, \textit{Duke Math.\ J.} \textbf{101} (2000),
359--426; \texttt{math/9908171}.

\bibitem{bib:KhR1}
M.\ Khovanov and L.\ Rozansky, Matrix factorizations and link homology
I, \texttt{math/0401268}.

\bibitem{bib:KhR2}
M.\ Khovanov and L.\ Rozansky, Matrix factorizations and link homology
II, \textit{Geom.\ Topol.}, to appear; \texttt{math/0505056}.

\bibitem{bib:KRKauff}
M.\ Khovanov and L.\ Rozansky, Virtual crossings, convolutions and a
categorification of the $SO(2N)$ Kauffman polynomial,
\texttt{math/0701333}.




\bibitem{bib:Mor}
H.\ R.\ Morton, Seifert circles and knot polynomials,
\textit{Math.\ Proc.\ Cambridge Philos.\ Soc.} \textbf{99} (1986),
no.\ 1, 107--109.


\bibitem{bib:NgKho}
L.\ Ng, A Legendrian Thurston--Bennequin bound from Khovanov
homology, \textit{Algebr.\ Geom.\ Topol.} \textbf{5} (2005),
1637--1653; \texttt{math/0508649}.

\bibitem{bib:Ngarc}
L.\ Ng, On arc index and maximal Thurston--Bennequin number,
\texttt{math/0612356}.




\bibitem{bib:OSz}
P.\ Ozsv\'ath and Z.\ Szab\'o, Knot Floer homology and the four-ball
genus, \textit{Geom.\ Topol.} \textbf{1} (2001), 427--434;
\texttt{math/0301149}.

\bibitem{bib:Plam1}
O.\ Plamenevskaya, Bounds for the Thurston--Bennequin number from
Floer homology, \textit{Algebr.\ Geom.\ Topol.} \textbf{4} (2004),
399--406; \texttt{math/0311090}.

\bibitem{bib:Plam2}
O.\ Plamenevskaya, Transverse knots and Khovanov homology,
\textit{Math.\ Res.\ Lett.} \textbf{13} (2006), no.\ 4, 571--586;
\texttt{math/0412184}.

\bibitem{bib:Ras}
J.\ Rasmussen, Khovanov homology and the slice genus, \texttt{math/0402131}.

\bibitem{bib:RasKR}
J.\ Rasmussen, Some differentials on Khovanov--Rozansky homology,
\texttt{math/0607544}.

\bibitem{bib:Rud1}
L.\ Rudolph, A congruence between link polynomials,
\textit{Math.\ Proc.\ Cambridge Philos.\ Soc.} \textbf{107}
(1990), 319--327.

\bibitem{bib:Rud2}
L.\ Rudolph, Quasipositivity as an obstruction to sliceness,
\textit{Bull.\ Amer.\ Math.\ Soc.\ (N.S.)} \textbf{29} (1993), 51--59.

\bibitem{bib:Ru}
D.\ Rutherford, Thurston--Bennequin number, Kauffman
polynomial, and ruling invariants of a Legendrian link: the Fuchs
conjecture and beyond, \textit{Int.\ Math.\ Res.\ Not.} \textbf{2006},
Art.\ ID 78591;
\texttt{math/0511097}.

\bibitem{bib:Shu}
A.\ Shumakovitch, Rasmussen invariant, slice-Bennequin inequality, and
sliceness of knots, \texttt{math/0411643}.


\bibitem{bib:Tab}
S.\ Tabachnikov, Estimates for the Bennequin number of Legendrian
links from state models for knot polynomials, \textit{Math.\ Res.\
  Lett.} \textbf{4} (1997), no.\ 1, 143--156.

\bibitem{bib:Tan}
T.\ Tanaka, Maximal Bennequin numbers and Kauffman polynomials of
positive links, \textit{Proc.\ Amer.\ Math.\ Soc.} \textbf{127}
(1999), no.\ 11, 3427--3432.

\bibitem{bib:Wu}
H.\ Wu, Braids, transversal links and the Khovanov--Rozansky homology,
\texttt{math/0508064}.

\bibitem{bib:Wu2}
H.\ Wu, Legendrian links and the spanning tree model for Khovanov
homology, \textit{Algebr.\ Geom.\ Topol.} \textbf{6} (2006),
1745--1757; \texttt{math/0605630}.

\bibitem{bib:Wu3}
H.\ Wu, The Khovanov--Rozansky cohomology and Bennequin inequalities,
\texttt{math/0703210}.

\end{thebibliography}
\end{document}